\documentclass{article}

\usepackage[all]{xy}
\usepackage{mathrsfs}
\usepackage{amssymb}
\usepackage{amsthm}
\usepackage{amsmath}
\usepackage{amsfonts}
\usepackage{authblk}

\usepackage{url}

\newtheorem{thm}{Theorem}

\allowdisplaybreaks

\author{Camilo Sanabria Malagón\thanks{c.sanabria135@uniandes.edu.co}}
\affil{Departamento de Matemáticas\\ Universidad de los Andes\\ Bogotá, Colombia}

\title{A note on the inverse problem for finite differential Galois groups}

\begin{document}

\date{}

\maketitle

\begin{abstract}
  In this paper we revisit the following inverse problem: given a curve invariant under an irreducible finite linear algebraic group, can we construct an ordinary linear differential equation whose Schwarz map parametrizes it? We present an algorithmic solution to this problem under the assumption that we are given the function field of the quotient curve. The result provides a generalization and an efficient implementation of the solution to the inverse problem exposed by M. van der Put, C. Sanabria and J.Top \cite{VANDERPUT2020}. As an application, we show that there is no hypergeometric equation with differential Galois group isomorphic to $H_{72}^{SL_3}$, thus completing Beuker and Heckman's answer \cite{BEUKERS1989} to the question of which irreducible finite subgroup of $SL_3(\mathbb{C})$ are the monodromy of a hypergeometric equation.
\end{abstract}

\section*{Introduction}

  We consider a linear ordinary differential equation (LODE) $L(y)=0$, where
  $$L(y)=y^{(n)}+a_{n-1}y^{(n-1)}+\ldots+a_1y'+a_0y$$
  and $a_0,a_1,\ldots,a_{n-1}\in\mathbb{C}(z)$. Let $U\subseteq\mathbb{C}$ be an open disk avoiding the singularities of the equation and let
  $$y_j:U\longrightarrow\mathbb{C},\quad j\in\{1,\ldots,n\}$$
  be a basis of solutions. The multivalued analytic extension of
  \begin{eqnarray*}
  U & \longrightarrow & \mathbb{P}^{n-1}(\mathbb{C})\\
  z & \longmapsto & \Big[y_1(z):\ldots:y_n(z)\Big]
  \end{eqnarray*}
  is the Schwarz map associated to our basis of solutions and its image is an analytic curve in $\mathbb{P}^{n-1}(\mathbb{C})$.  The monodromy group of $L(y)=0$ defines a group of automorphisms of this curve.

  If all the solutions are algebraic over $\mathbb{C}(z)$, then the monodromy group is finite and the corresponding analytic curve is algebraic. The inverse problem can be stated in the following terms. Let $G\subset GL_n(\mathbb{C})$ be a finite linear algebraic group and $PG$ its projection in $PGL_n(\mathbb{C})\simeq GL_n(\mathbb{C})/Z$. Given an algebraic curve $C\subseteq \mathbb{P}^{n-1}(\mathbb{C})$, invariant under $PG$, can we construct a linear ordinary differential equation such that, for some basis of solutions, the image of its Schwarz map is the given curve and its monodromy group contains $G$? In the case where we have an affirmative answer, we say that the Schwarz map parametrizes $C$.

  One of the first results in this direction was given by G.H. Halphen \cite{HALPHEN1884} and A. Hurwitz \cite{HURWITZ1886} who independently solved the inverse problem in the case of the Klein quartic and its automorphism group $G_{168}$. Subsequent results arose from considering the same group but other invariant curves. For example, R. Fricke \cite{FRICKE1893} described a pencil of projective curves of degree $12$ invariant under $G_{168}$ and for each curve in this pencil M. Kato \cite{KATO2004} constructed an ordinary linear differential equation with its Schwarz map parametrizing it. Furthermore, M. Kato\cite{KATO2006} showed that among the hypergeometric equations with finite monodromy described by F. Beukers and G. Heckmann in \cite{BEUKERS1989}, there is one that solves the inverse problem for the Hessian of the Klein quartic.

  More generally, curves invariant under other finite linear algebraic groups have been studied. The Hesse pencil \cite{HESSE1844I, HESSE1844II}, a family of elliptic curves invariant under the normal abelian subgroup of order $9$ of the Hessian group is one example. This pencil has regained much attention recently in different areas \cite{TOP2018, ARTEBANI2009, SMART2001, ZASLOW2005}. For each curve in it C. Sanabria \cite{SANABRIA2017} constructed a linear ordinary differential equation with Schwarz map parametrizing it. Also, L. Lachtin \cite{LACHTIN1898} considered curves invariant under the Valentiner group, the perfect triple cover of $A_6$, and obtained a hypergeometric equation that parametrizes one of these curves.

  An important result closely related to this inverse problem is a classical theorem of Klein \cite{KLEIN1878,KLEIN1879} stating that, given a finite primitive group $G\subseteq SL_2(\mathbb{C})$, there exists a hypergeometric equation such that every second order LODE whose projective monodromy group is isomorphic to $PG$ is projectively equivalent to the pullback by a rational map of this hypergeometric equation. Indeed, for each finite primitive group $G\subseteq SL_2(\mathbb{C})$, there is a Gauss hypergeometric equation that parametrizes $\mathbb{P}^1(\mathbb{C})$ and its monodromy contains $G$. Naturally, these hypergeometric equations are in correspondence with Schwarz spherical triangles. This classical result has been generalized to LODEs of order $3$ by M. Berkenbosch \cite{BERKENBOSCH2006}, and to arbitrary order by C. Sanabria \cite{SANABRIA2014}. 

  A generalization of this inverse problem is treated by Y. Haraoka and M. Kato in \cite{HARAOKA2010}, where they consider projective spaces $\mathbb{P}^n(\mathbb{C})$ and reflection groups instead of curves and finite linear algebraic groups. Their solution provides another generalization of Klein's, one that applies to Fuchsian systems of rank $n$ over $\mathbb{P}^n(\mathbb{C})$ with monodromy isomorphic to a reflection group.

  A related inverse problem can be stated as follows: given a finite linear algebraic group $G\subseteq GL_n(\mathbb{C})$, can we construct a linear ordinary differential equation such that the monodromy group is isomorphic to $G$? This inverse problem was studied by F. Beukers and G. Heckman in \cite{BEUKERS1989} for reflection groups and by M. van der Put and F. Ulmer in \cite{VANDERPUT2000} for equations of order $3$. Their approaches obtain only equations over $\mathbb{P}^1(\mathbb{C})$ with three singularities, and couldn't determine whether there is a hypergeometric equation with differential Galois group $H_{72}^{SL_3}$.

  In \cite{VANDERPUT2020}, M. van der Put, C. Sanabria and J. Top presented a solution to the inverse problem in the case when the normalization of the quotient of the curve by the group of $PG$ has genus $0$. More precisely, given any algebraic curve $C\subseteq \mathbb{P}^{n-1}(\mathbb{C})$ invariant under a finite linear algebraic group $G\subseteq GL_n(\mathbb{C})$ together with a generator of $\mathbb{C}(C/G)$, they presented a procedure that constructs a linear ordinary differential equation whose Schwarz map parametrizes $C$.  Furthermore, they showed that any irreducible linear ordinary differential equation with rational coefficients and only algebraic solutions can be obtained as a result of this procedure.

  In this paper, we revisit the solution in \cite{VANDERPUT2020} to the inverse problem but from a geometric perspective. In this way, we obtain a generalization where we drop the condition on the genus of the normalization of the quotient curve. Moreover, we turn the solution into an efficient algorithm. Finally, as an application of this alternative viewpoint, we prove that there is no hypergeometric equation with differential Galois group isomorphic to the elusive group $H_{72}^{SL_3}$, thus completing Beuker and G. Heckman answer \cite{BEUKERS1989} to the question of which irreducible finite subgroup of $SL_3(\mathbb{C})$ are the monodromy of a hypergeometric equation.

\section{Preliminary}

\subsection{Space of orbits}\label{orbitspace}

  We recall some elements of the space of orbits of actions of finite linear groups on
  vector spaces. A more detailed exposition is found in \cite{BRION2010}.

  Let $G\subseteq GL_n(\mathbb{C})$ be a finite linear algebraic group. Let $R=\mathbb{C}[X_1,\ldots,X_n]$ be the coordinate ring of $\mathbb{C}^n$. We consider the natural left action of $G$ on $R$ by $\mathbb{C}$-automorphisms
  \[
  g: X_j\mapsto \sum_{i=1}^n X_ig_{ij},
  \]
  for $g=(g_{ij})_{i,j=1}^n\in G$.
  This $G$-action on $R$ induces a right $G$-action on $\mathbb{C}^n$ defined by
  \[
  g: (x_1,\ldots,x_n)\mapsto \Big(\sum_{i=1}^nx_ig_{i1},\ldots,\sum_{i=1}^nx_ig_{in}\Big).
  \]
  Since $G$ is finite, it is reductive, thus the $G$-invariant polynomials in $R$ separate the $G$-orbits in $\mathbb{C}^n$. Therefore, the coordinate ring of the orbit space $\mathbb{C}^n/G$ is the finitely generated $G$-invariant subring $R^G\subseteq R$.

  Let $F_1,\ldots,F_N$ be homogeneous generators of $R^G$. We can embed $\mathbb{C}^n/G$ into $\mathbb{C}^N$ through the algebraic map
  \[
  \mathbf{x}\cdot G \mapsto \Big(F_1(\mathbf{x}),\ldots,F_N(\mathbf{x})\Big),
  \]
  where $\mathbf{x}=(x_1,\ldots,x_n)$. We identify the orbit space $\mathbb{C}^n/G$ with its embedding in $\mathbb{C}^N$. Since $G$ is finite, the orbit space has dimension $n$ and therefore there exists an open dense subset of $\mathbb{C}^n$ where the derivative of the quotient map
  \[
  \mathbf{x}\mapsto \Big(F_1(\mathbf{x}),\ldots,F_N(\mathbf{x})\Big)
  \]
  is non-singular.

  Similarly,  we get natural $G$-actions on $\mathbb{P}^{n-1}(\mathbb{C})$ and on its homogeneous coordinate ring $R$, and we denote $\mathbb{P}^{n-1}(\mathbb{C})/G$ the orbit space. We denote by $\Pi: \mathbb{P}^{n-1}(\mathbb{C})\rightarrow \mathbb{P}^{n-1}(\mathbb{C})/G$ the quotient map. Let $\Lambda\in\mathbb{Z}_{>0}$ be such that $R^G_\Lambda$, the homogeneous elements of $R^G$ of degree $\Lambda$, form a homogeneous coordinate system for $\mathbb{P}^{n-1}(\mathbb{C})/G$, i.e. $$\mathbb{P}(R^G_\Lambda)\simeq\mathbb{P}^{n-1}(\mathbb{C})/G.$$
  Let $\Phi_1,\ldots,\Phi_M$ be a basis of $R^G_\Lambda$. We can embed $\mathbb{P}^{n-1}(\mathbb{C})/G$ into $\mathbb{P}^{M-1}(\mathbb{C})$ with the map
  \[
  [\mathbf{x}]\cdot G\mapsto \Big[\Phi_1(\mathbf{x}):\ldots:\Phi_M(\mathbf{x})\Big],
  \]
  where $[\mathbf{x}]=[x_1:\ldots:x_n]$. As before, we identify the orbit space $\mathbb{P}^{n-1}(\mathbb{C})/G$ with its embedding in $\mathbb{P}^{M-1}(\mathbb{C})$. Again, there exists an open dense subset of $\mathbb{P}^{n-1}(\mathbb{C})$ where the derivative of the quotient map
  \[
  [\mathbf{x}]\mapsto \Big[\Phi_1(\mathbf{x}):\ldots:\Phi_M(\mathbf{x})\Big]
  \]
  is non-singular.

\subsection{Schwarz maps}

  We now introduce Schwarz maps. A famous reference on the subject is \cite{YOSHIDA1997}

  Let $C_0$ be a compact Riemann surface and $K=\mathbb{C}(C_0)$ the field of meromorphic functions over $C_0$. Let $\delta:K\rightarrow K$ be any non-trivial derivation. Note that $\delta$ can be uniquely extended to the sheaf of meromorphic functions over open sets of $C_0$, and in particular to any algebraic extension of $K$.

  Let $L(y)=0$ be an linear ordinary  differential equation of order $n$ where
  \[
  L(y)=\delta^n(y)+a_{n-1}\delta^{n-1}(y)+\ldots+a_1\delta(y)+a_0y,
  \]
  and $a_0,a_1,\ldots,a_{n-1}\in K$. Let $S_0\subset C_0$ be the set of singular points of $L(y)$. Given a non-singular point $p\in C_0$ together with a basis of solutions $\mathbf{y}=(y_1,\ldots,y_n)$ defined over a neighborhood $U\subseteq C_0$ of $p$, we define the Schwarz map as the analytic extension of 
  \begin{eqnarray*}
  [\mathbf{y}]: U & \longrightarrow & \mathbb{P}^{n-1}(\mathbb{C})\\
  z & \longmapsto & \Big[y_1(z):\ldots:y_n(z)\Big].
  \end{eqnarray*}

  The monodromy of the Schwarz map is the projection on $PGL_n(\mathbb{C})$ of the monodromy group of the basis of solutions $(y_1,\ldots,y_n)$ . We will denote it by $PG$ and call it the projective monodromy of $L(y)=0$. Let us assume $PG$ is finite. Then, post-composing the Schwarz map with the quotient by the action of $PG$ on  $\mathbb{P}^{n-1}(\mathbb{C})$ we obtain a (single-valued) rational map $$\psi:C_0\dashrightarrow\mathbb{P}^{n-1}(\mathbb{C})/PG,$$ defined over $C_0\setminus S_0$, which we call the quotient Schwarz map.

\subsection{Algebraic Picard-Vessiot theory}

  We briefly recall the results from Picard-Vessiot theory that we will use in this paper. An extensive exposition of this subject is given in \cite{VANDERPUT2003}.

  We assume that $L(y)=0$ is irreducible and that all its solutions are algebraic over $K$. Under these assumptions, if $(y_1,\ldots,y_n)$ is a basis of solutions, then a Picard-Vessiot extension of $K$ for $L(y)=0$ is the field
  $$E=K[y_1,\ldots,y_n].$$
  Let's denote by $I$ the kernel of the evaluative $K$-morphism
  \begin{eqnarray*}
  \Phi: K[Y_1,\ldots,Y_n] & \longrightarrow & K\\
  Y_j & \longmapsto & y_j.
  \end{eqnarray*}
  An isomorphic representation of the differential Galois group of $L(y)=0$ is the finite group $G\subset GL_n(\mathbb{C})$ composed by the elements $g=(g_{ij})_{i,j=1}^n$ that send the ideal $I$ into itself under the $K$-automorphism
  \begin{eqnarray*}
  g: K[Y_1,\ldots,Y_n]  & \longrightarrow & K[Y_1,\ldots,Y_n]\\
  Y_j & \longmapsto & \sum_{i=1}^nY_ig_{ij}.
  \end{eqnarray*}
  By the Galois correspondence, if $P\in K[Y_1,\ldots,Y_n]$ is $G$-invariant then $\Phi(P)\in K$. Since $L(x)=0$ is irreducible, then $G$ is reductive and we have the following result.

  \begin{thm}\label{compointthm}[Compoint's theorem for the algebraic case \cite[Theorem 2]{SANABRIA2017}\cite[Corollary 4.2]{VANDERPUT2020}] Let $L(y)=0$ be an irreducible linear ordinary differential equation. If the differential Galois group of $L(y)=0$ is finite, then the ideal $I$ is generated by the $G$-invariants contained in it. In particular, if $P_1,\ldots,P_N\in\mathbb{C}[Y_1,\ldots,Y_n]$ is a set of generators of the $G$-invariant subring $\mathbb{C}[Y_1,\ldots,Y_n]^G$, then
  \begin{equation*}\label{Igen}
  I=\langle P_1-f_1,\ldots, P_N-f_N\rangle,
  \end{equation*}
  where $f_i=\Phi(P_i)$, $i=1,\ldots,N$.
  \end{thm}

\section{Schwarz maps for invariant curves}

  Let $C\subseteq\mathbb{P}^{n-1}(\mathbb{C})$ be an algebraic $G$-invariant irreducible curve. Let $C_0\rightarrow C/G$ be a normalization of the quotient curve. The following theorem shows that we can find a linear ordinary differential equation, with coefficients in an abelian extension of $\mathbb{C}(C_0)$, whose Schwarz map parametrizes $C$. The result provides a generalization of \cite[Theorem 6]{SANABRIA2017} and \cite[Procedure 4.3]{VANDERPUT2020}.



\begin{figure}
    {\center
    \hfill
    \xymatrix{
      & & C\ \ar[ddd]\ar@{^{(}->}[r]_{\tilde{\imath}} & \mathbb{P}^{n-1}(\mathbb{C})\ar[ddd]^{\Pi}\\
      & & \\
    C_1\ar@{..>}[rruu]^{[\mathbf{y}]}\ar[d]_\pi  &  & \\
    C_0\ar@{-->}[rr]_\sigma^\simeq & &  C/G\ \ar@{^{(}->}[r]_\imath & \mathbb{P}^{n-1}(\mathbb{C})/G
    }
    \hfill
    }
    \caption{Diagram for Theorem \ref{theo}\label{theodia}}
  \end{figure}
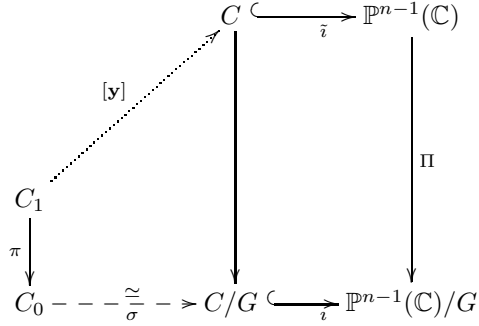

  \begin{thm}\label{theo}
    Let $G$ be a finite algebraic subgroup of $GL_n(\mathbb{C})$ and let $C\subseteq\mathbb{P}^{n-1}(\mathbb{C})$ be an algebraic $G$-invariant irreducible curve not contained in a projective hyperplane. Let $\sigma:C_0\rightarrow C/G$ be a normalization. Let $K=\mathbb{C}(C_0)$ be the field of meromorphic functions of $C_0$. Then there exist (see Figure \ref{theodia})
    \begin{itemize}
    \item[i)] a branched cover $\pi:C_1\rightarrow C_0$ where $C_1$ is a compact Riemann surface and $K\subseteq E=\mathbb{C}(C_1)$ is an abelian extension;
    \item[ii)] a linear differential equation $L(y)=0$ of order $n$ with coefficients in $E$ whose Schwarz map parametrizes $C$ and the quotient Schwarz map is $\psi=\sigma\circ\pi$.
    \end{itemize}
  \end{thm}

  \begin{proof}
    We recall our notation from section \ref{orbitspace}. The invariant ring $R^G=\mathbb{C}[X_1,\ldots,X_n]^G$, is generated by the homogeneous polynomials $F_1,\ldots,F_N$. The polynomials $\Phi_1,\ldots,\Phi_M$ form a basis of the homogeneous coordinate system $R^G_\Lambda$ for $\mathbb{P}^{n-1}(\mathbb{C})/G$.

    For $j=1,\ldots,M$, let $r_{j,1},\ldots,r_{j,N}\in\mathbb{Z}_{\ge 0}$ be such that $r_{j,1}\deg\left(F_1\right)+\ldots+r_{j,N}\deg\left(F_N\right)=\Lambda$ and
    $$\Phi_j=\prod_{i=1}^N F_i^{r_{j,i}}.$$

    Let $\tilde{\imath}: C\rightarrow \mathbb{P}^{n-1}(\mathbb{C})$ denote the inclusion map. Likewise, let $\imath: C/G\rightarrow \mathbb{P}^{n-1}(\mathbb{C})/G$ be the inclusion and let $d\in\{1,\ldots,M\}$ be such that $C/G\not\subseteq\{\Phi_d=0\}$. Then there exist $\phi_1,\ldots, \phi_N\in K$ such that for $i=1,\ldots,N$
    $$(\imath\circ\sigma)^*\left(\dfrac{\Phi_i}{\Phi_d}\right)=\dfrac{\phi_i}{\phi_d}.$$
    Let $f_1,\ldots,f_N\in\overline{K}$ be such that
    $$\phi_j=\prod_{i=1}^N f_i^{r_{j,i}}.$$
    Then, the field $E=K(f_1,\ldots,f_n)$ is an abelian extension of $K$. Let $C_1$ be a compact Riemann surface with field of meromorphic functions $\mathbb{C}(C_1)\simeq E$ and let $\pi:C_1\rightarrow C_0$ be the map defined by the inclusion $K\subseteq E$.

    Let $U\subseteq C/G$ be an open set not containing any ramification point of $C\rightarrow C/G$, nor any singularity of $C/G$. Furthermore, we may take $U$ such that, if $U_1=(\sigma\circ\pi)^{-1}(U)$, then $\sigma\circ\pi\upharpoonright_{U_1}$ is a bijection. Let $y_1,\ldots y_n$ be meromorphic functions on $U_1$ such that $[\mathbf{y}]: z\mapsto [y_1(z):\ldots:y_n(z)]$ is a lifting of $\sigma\circ\pi$, i.e. $\Pi\upharpoonright_C\circ[\mathbf{y}]=\sigma\circ\pi$. So, for $j=1,\ldots,M$, we have that
    \begin{align*}
      \dfrac{\Phi_j}{\Phi_d} \circ \tilde{\imath} \circ[\mathbf{y}] &
        = \dfrac{\Phi_j}{\Phi_d} \circ \Pi \circ \tilde{\imath}\circ [\mathbf{y}]
        = \dfrac{\Phi_j}{\Phi_d} \circ \imath \circ \Pi\upharpoonright_C \circ [\mathbf{y}] \\
      & = \dfrac{\Phi_j}{\Phi_d} \circ \imath \circ \sigma \circ \pi
        = \pi^*\circ (\imath \circ \sigma)^*\left(\dfrac{\Phi_j}{\Phi_d}\right) = \pi^* (\dfrac{\phi_j}{\phi_d}) = \dfrac{\phi_j}{\phi_d},
    \end{align*} 
    and therefore
    $$\dfrac{\Phi_i}{\Phi_d}(y_1,\ldots,y_n)=\dfrac{\phi_i}{\phi_d}.$$
    Replacing $(y_1,\ldots,y_n)\leftarrow \lambda(y_1,\ldots,y_n)$ if necessary, where $\lambda$ is a radical of an element in $K$, we may assume $\Phi_j(y_1,\ldots,y_n)=\phi_j$ for $j=1,\ldots,M$, and furthermore, for $i=1,\ldots,N$
    \begin{equation}\label{Fifi}
      F_i(y_1,\ldots,y_n)=f_i.
    \end{equation}
    Let $\delta:K\rightarrow K$ be a non-trivial derivation. Differentiating $\eqref{Fifi}$, we obtain
    \[
    \frac{\partial F_i}{\partial X_1}(\mathbf{y})\delta(y_1)+\ldots+\frac{\partial F_i}{\partial X_n}(\mathbf{y})\delta(y_n)=  \delta(f_i),
    \]
    where $\mathbf{y}=(y_1,\ldots,y_n)$. Since the derivative of $\mathbf{x}\mapsto \Big(F_1(\mathbf{x}),\ldots,F_N(\mathbf{x})\Big)$ is non-singular in a dense open set of $\mathbb{C}^n$, after rearranging indices if necessary, we have that
    $$\frac{\partial(F_1,\ldots,F_n)}{\partial(X_1,\ldots,X_n)}(\mathbf{x})=\left[\begin{array}{ccc}
    \partial F_1/\partial X_1 & \ldots & \partial F_1/\partial X_n\\
    \vdots & \ddots & \vdots\\
    \partial F_n/\partial X_1 & \ldots & \partial F_n/\partial X_n
    \end{array}\right](\mathbf{x})$$
    is invertible in that set. Whenever $\det\Big(\frac{\partial(F_1,\ldots,F_n)}{\partial(X_1,\ldots,X_n)}\Big)(\mathbf{y}(z))\ne 0$ for $z\in U_1$, we obtain a system of first order ordinary differential equations characterizing the functions $y_1,\ldots ,y_n$
    \begin{equation}\label{ode}
    \delta\left[\begin{array}{c}
    y_1 \\ \vdots \\ y_n
    \end{array}\right] =
    \left[\frac{\partial(F_1,\ldots,F_n)}{\partial(X_1,\ldots,X_n)}(\mathbf{y})\right]^{-1}\ \delta\left[\begin{array}{c}
    f_1 \\ \vdots \\ f_n
    \end{array}\right].
    \end{equation}
    Differentiating this system $n-1$ times, we obtain expressions for $\delta^k\mathbf{y}$, $k=2,\ldots,n$,
    in terms of the derivatives of $f_1,\ldots, f_n$ and the partial derivatives of $F_1,\ldots,F_n$. The system \eqref{ode} is $G$-invariant and therefore so are the linear dependence relations
    $$c_0\mathbf{y} + c_1\delta\mathbf{y} + \quad \ldots \quad + c_{n-1}\delta^{n-1}\mathbf{y} + c_n\delta^n\mathbf{y} = 0,$$
    hence the $n+1$ coefficients $c_0,\ldots,c_n$ are in the field
    $$\left(F_i(\mathbf{y}),\delta^k f_j\Big|\ i,j,k=1,\ldots,n\right)\subseteq\mathbb{C}\left(\delta^kf_i\Big|\ i=1,\ldots,N,\ k=0,1,\ldots,n\right)\subseteq E.$$
    We have that $c_n\ne 0$ because $C$ is not contained in a projective hyperplane. So, the linear differential equation $$\delta^n(y)+\frac{c_{n-1}}{c_n}\delta^{n-1}(y)+\ldots+\frac{c_1}{c_n}\delta(y)+\frac{c_0}{c_n}y=0$$ has coefficients in $E$ and $(y_1,\ldots,y_n)$ as a basis of solutions. By construction, the analytic extension of $z\mapsto [\mathbf{y}]=[y_1(z):\ldots:y_n(z)]$ is a Schwarz map that parametrizes $C$ and the quotient Schwarz map is $\sigma\circ\pi$.
  \end{proof}

\subsection*{The algorithm}

  We present an algorithm that computes the equation $L(y)=0$ from Theorem \ref{theo}. We assume that $K\subseteq\overline{\mathbb{Q}}(z)$ is a field such that we can compute Gröbner basis on $K[X_1,\ldots,X_n]$. Furthermore, we also assume that the homogeneous generators $F_1,\ldots,F_N$ of $\mathbb{C}[X_1,\ldots,X_n]^G$ have coefficients in $\mathbb{K}=K\cap\overline{\mathbb{Q}}$ and that they are indexed so that $\det\Big(\frac{\partial(F_1,\ldots,F_n)}{\partial(X_1,\ldots,X_n)}\Big)\ne 0$, as in the case where $F_1,\ldots,F_n$ form a homogeneous system of parameters of $\mathbb{C}[X_1,\ldots,X_n]^G$. Details on Gröbner bases and homogeneous system of parameters can be found in \cite{DERKSEN2015}.

  We take $\delta=d/dz$. 


  \noindent\textbf{Algorithm}\newline
  \textbf{Input}: $f_1,\ldots,f_N\in K$ such that 
  \begin{itemize}
    \item[i)]$I=\langle F_1-f_1,\ldots, F_N-f_N\rangle$ is a maximal ideal of $K[X_1,\ldots,X_n]$, and
    \item[ii)] the projective curve $C\subseteq\mathbb{P}^{n-1}(\mathbb{C})$ defined by the homogeneous elements in $I\cap\mathbb{K}[X_1,\ldots,X_n]$ is not contained in a projective hyperplane.
  \end{itemize}
  \textbf{Output}: $(a_{n-1},\ldots,a_1,a_0)\in K^n$ such that the linear ordinary differential equation $\delta^n(y)+a_{n-1}\delta^{n-1}(y)+\ldots+a_1\delta(y)+a_0y=0$ has a fundamental system of solutions whose associated Schwarz map parametrizes $C$.
  \begin{enumerate}
    \item Fix an admissible ordering on the monomials on $X_1,\ldots,X_n$ and construct a Gr\"obner basis $B$ of $\langle F_1-f_1,\ldots, F_N-f_N\rangle$.
    \item For $i=1,\ldots,n$, set $Y_{0,i}:=X_i$. Denote $\mathbf{Y}_0=(Y_{0,1},\ldots,Y_{0,n})^\intercal$.
    \item Let $\mathbf{J}^{-1}$ be the remainder by $B$ of an inverse $\bmod I$ of the jacobian matrix $\partial(F_1,\ldots,F_n)/\partial(X_1,\ldots,X_n)$.
    \item Let $Y_{1,1},\ldots,Y_{1,n}$ be such that $\mathbf{Y}_1$ is the remainder by $B$ of $\mathbf{J}^{-1}\mathbf{f}$ where $\mathbf{Y}_1=[Y_{1,1},\ldots,Y_{1,n}]^\intercal$ and $\mathbf{f}=[f_1,\ldots,f_n]^\intercal$.
    \item Extend $\delta$ to a derivative in $K(X_1,\ldots,X_n)$ by defining $\delta X_i=0$ for $i=1,\ldots,n$.
    \item For $j=2,\ldots,n$ recursively let $Y_{j,1},\ldots,Y_{j,n}$ be such that, for $i=1,\ldots,n$, $Y_{j,i}$ is the remainder by $B$ of 
    $$\delta Y_{j-1,i}+\sum_{k=1}^n\dfrac{\partial {Y_{j-1,i}}}{X_k}Y_{1,k}.$$
    Denote $\mathbf{Y}_j=(Y_{j,1},\ldots,Y_{j,n})^\intercal$.
    \item Let $\mathbf{Y}^{-1}$ be the remainder by $B$ of an inverse $\bmod I$ of $[\mathbf{Y}_{n-1}|\ldots|\mathbf{Y}_1|\mathbf{Y}_0]$.
    \item \textbf{Return} the remainder $\mathbf{a}$ by $B$ of $-\mathbf{Y}^{-1}\mathbf{Y}_n$.
  \end{enumerate}

  We note that, as a corollary of Compoint's theorem, we have that every irreducible linear ordinary differential equation with coefficients in $\mathbb{Q}(z)$ and finite differential Galois group can be obtained from the algorithm.

\section{The group $H_{72}^{SL_3}$}

  The group $H_{216}^{SL_3}$ is the preimage in $SL_3(\mathbb{C})$ of the Hessian group $H_{216}\subseteq PSL_3(\mathbb{C})$.

  As in \cite{SINGER1993}, we defined the following matrices.
  $$
  S=\left[\begin{array}{ccc}
  1 & 0 & 0\\
  0 & \omega & 0\\
  0 & 0 & \omega^2
  \end{array}\right], \quad
  T=\left[\begin{array}{ccc}
  0 & 1 & 0\\
  0 & 0 & 1\\
  1 & 0  & 0
  \end{array}\right], \quad
  U=\left[\begin{array}{ccc}
  \varepsilon & 0 & 0\\
  0 & \varepsilon & 0\\
  0 & 0  & \varepsilon\omega
  \end{array}\right],
  $$
  and
  $$
  V=\rho\left[\begin{array}{ccc}
  1 & 1 & 1\\
  1 & \omega & \omega^2\\
  1 & \omega^2  & \omega
  \end{array}\right].
  $$
  Here, $\varepsilon$ is a primitive $9$th root of unity. Thus, $\varepsilon^6+\varepsilon^3+1=0$ and $\omega=-1-\varepsilon^3$ is a primitive $3$rd root of unity. The scaling factor $\rho=1/(\omega-\omega^2)$ is defined so that $\det(V)=1$. The generating matrices 

  The group $H_{216}^{SL_3}$ is generated by the matrices $S, T, U, V$. The Hessian group has a normal subgroup of index $3$, $H_{72}$. The preimage of this subgroup is $H_{72}^{SL_3}$ and it's generated by $S$, $T$, $V$, $UVU^{-1}$. The group $H_{72}$ has a normal subgroup of index $2$, $F_{36}$. The preimage of this subgroup is $F_{36}^{SL_3}$, and it's generated by $S$, $T$ and $V$.

\subsubsection*{Invariants of $H_{72}^{SL_3}$}

  As in \cite{ROTILLON1981}, we define $P =X_1X_2X_3$, $S = X_1^3+X_2^3+X_3^3$, $Q=X_1^3X_2^3+X_1^3X_3^3+X_2^3X_3^3$, $R = (X_1^3-X_2^3)(X_1^3-X_3^3)(X_2^3-X_3^3)$,
  $F_6 = S^2-12Q$, $\Phi_6 =S^2-18P^2-6PS$, and $F_{12} =S^4+216P^3S$.

  The ring of invariants of $H_{72}^{SL_3}$ is
  $$\mathbb{C}[X_1,X_2,X_3]^{H_{72}^{SL_3}}  =\mathbb{C}[F_6,R,F_{12},\Phi_6^2].$$
  As a ring, it's isomorphic to
  $$\mathbb{C}[Z_6,Z_9,Z_{12},X_{12}]/(T_{36}),$$
  where
  \begin{align*}
  T_{36}  & =(432Z_9^2+3Z_6Z_{12}-Z_6^3)^2-4(X_{12}^3-3Z_{12}X_{12}^2+3Z_{12}^2X_{12})\\
    & =(432Z_9^2+3Z_6Z_{12}-Z_6^3)^2-4X_{12}\big((\xi+1)Z_{12}-X_{12}\big)\big((\overline{\xi}+1)Z_{12}-X_{12}\big)
  \end{align*}
  and $\xi$ is a $6$th root of unity.

  The polynomials $F_6,R,F_{12}$ form a homogeneous system of parameters of $\mathbb{C}[X_1,X_2,X_3]^{H_{72}^{SL_3}}$ \cite[Théorème 4]{ROTILLON1981}. The jacobian of the map $(X_1,X_2,X_3) \mapsto (F_6,R,F_{12})$ is
  \begin{eqnarray*}
  \det\left(\frac{\partial(F_6,R,F_{12})}{\partial(X_1,X_2,X_3)}\right) & = & \left|\begin{array}{ccc}
  \partial F_6/\partial X_1 & \partial F_6/\partial X_2 & \partial F_6/\partial X_3\\
  \partial R/\partial X_1 & \partial R/\partial X_2 & \partial R/\partial X_3\\
  \partial F_{12}/\partial X_1 & \partial F_{12}/\partial X_2 & \partial F_{12}/\partial X_3\\
  \end{array}\right|\\
  & = & 18(F_{12}-\Phi_6^2)^2
  \end{eqnarray*}

  The ring of invariants of $F_{36}^{SL_3}$ is
  $$\mathbb{C}[X_1,X_2,X_3]^{F_{36}^{SL_3}}  =\mathbb{C}[F_6,\Phi_6,R,F_{12},\Psi_{12}],$$
  where $\Psi_{12}=PS^3+3P^2S^2-18P^3S$. As a ring, it's isomorphic to
  $$\mathbb{C}[Z_6,X_6,Z_9,Z_{12},X_{12}]/(T_{18}, T_{24}),$$
  where
  \begin{align*}
  T_{18}  & =432 Z_9^2 - Z_6^3 + 3 Z_6 Z_12 - 2 X_6^3 - 36 X_6 X_{12},\\
  T_{24}  & =Z_{12}^2 - X_6^2 (Z_{12} + 12 X_{12}) + 12 X_{12}^2
  \end{align*}

  For future reference, we remark that $\Phi_6$ factors into semi-invariants of $F_{36}^{SL_3}$ as $\Phi_6 = (S + \xi^2P)(S + \xi^{-2}P)$. The factors are semi-invariants of order $4$, and the $\mathbb{C}$-vector space generated by them is irreducible under the action of $H_{72}^{SL_3}$.

\subsection*{Differential equations for $H_{72}^{SL_3}$}

  F. Beuker and G. Heckman in \cite{BEUKERS1989} and by M. van der Put and F. Ulmer in \cite{VANDERPUT2000} obtained different families of hypergeometric equations of degree $3$ with finite differential Galois group. They have hypergeometric equations for each finite primitive subgroup of $SL_3(\mathbb{C})$ but for $H_{72}^{SL_3}$. We now show that there isn't any.

  From Theorem \ref{theo}, to an LODE with differential Galois group isomorphic to $H_{72}^{SL_3}$ corresponds a $\mathbb{C}$-algebra homomorphism $\psi^*: \mathbb{C}[Z_6,Z_9,Z_{12},X_{12}]/(T_{36}) \rightarrow \mathbb{C}[z]$ induced by the quotient Schwarz map. We denote by $f_6$, $r_9$, $f_{12}$ and $\phi_6^2$ the images of $Z_6$, $Z_9$, $Z_{12}$ and $X_{12}$ respectively.

  Let's first see that that $\phi_6^2\ne 0$. Indeed, if $\phi_6^2= 0$, we have $\Phi_6(\mathbf{y})=0$ for some basis of solutions $\mathbf{y}=(y_1,y_2,y_3)$ of the LODE. Whence $(S + \xi^2P)(\mathbf{y})=0$ or $(S + \xi^{-2}P)(\mathbf{y})=0$, and so $S + \xi^2P$ and $S + \xi^{-2}P$ are semi-invariants of the Galois group, and therefore, the Galois group would be isomorphic to a subgroup of $F_{36}^{SL_3}$ instead of $H_{72}^{SL_3}$.

  Let's now see that $\phi_{6}^2$ must have at least one zero of odd multiplicity. For, if not, we would have $\phi_6\in\mathbb{C}[z]$ and, since $\phi_6^2\ne 0$, from the syzygy $T_{18}$, we have $\psi_{12}=(432 r_9^2 - f_6^3 + 3 f_6 f_{12} - 2 \phi_6^3)/(36 \phi_6)\in\mathbb{C}(z)$, therefore $\Phi_6$ and $\Psi_{12}$ would be invariants and, again, the Galois group would be isomorphic to a subgroup of $F_{36}^{SL_3}$ instead of $H_{72}^{SL_3}$.

  Now, since $\psi^*(T_{36})=0$, we have
  \begin{align*}
    (432r_9^2+3f_6f_{12}-f_6^3)^2 & = 4(\phi_{6}^6-3f_{12}\phi_{6}^4+3f_{12}^2\phi_{6}^2)\\
      & = 4\phi_{6}^2\big(\phi_{6}^2-(\xi+1)f_{12}\big)\big(\phi_{6}^2-(\overline{\xi}+1)f_{12}\big).
  \end{align*}
  The right-hand side is a square in $\mathbb{C}[z]$, and, in particular, each of its zeroes has even multiplicity. So, if $z_0\in\mathbb{C}$ is a zero of odd multiplicity of one of the three factors $\phi_{6}^2$, $(\xi+1)f_{12}-\phi_{6}^2$ or $(\overline{\xi}+1)f_{12}-\phi_{6}^2$, then it must be a zero of odd multiplicity of another of these three factor but not of both. But since any of the three factors is a linear combination over $\mathbb{C}$ of the two others, if $z_0$ is a zero of two factors, it's a zero of the third factor too. Thus, any zero of odd multiplicity of one of the three factors, since it's a zero of odd multiplicity of two factors, then it's a zero of even multiplicity of the third factor. Similarly, a zero of two factors is a zero of $\phi_6^2$ and $f_{12}$ because these two values are linear combinations over $\mathbb{C}$ of any of the two factors.

  Let $z_0$ be a zero of odd multiplicity of $\phi_6^2$. Assume first that $z_0$ is a zero of odd multiplicity of the factor $\phi_{6}^2-(\xi+1)f_{12}$, and so a zero of even multiplicity of the factor $\phi_{6}^2-(\overline{\xi}+1)f_{12}$. Since $\phi_{6}^2-(\xi+1)f_{12}$ and $\phi_{6}^2-(\overline{\xi}+1)f_{12}$ are conjugates, $\overline{z_0}$ is a zero of odd multiplicity of $\phi_{6}^2-(\overline{\xi}+1)f_{12}$ and a zero of even multiplicity of $\phi_{6}^2-(\xi+1)f_{12}$. Therefore, $z_0\ne \overline{z_0}$, and in particular $z_0\ne 0, 1$. Since $z_0$ and $\overline{z_0}$ are zeroes of two of the factors, $z_0$ is a zero of $\phi_6^2$ and of $f_{12}$, and so of the jacobian $18(f_{12}-\phi_6^2)^2$. Hence, $z_0$ is a singularity of the equation different from $0$ and $1$, and so the equation is not hypergeometric. And similarly, if we assume that $z_0$ is a zero of odd multiplicity of $\phi_{6}^2-(\overline{\xi}+1)f_{12}$. As a corollary, we obtain that every third order LODE with differential Galois group $H_{72}^{SL_3}$ and rational coefficients has non-real singularities. Two such equations can be found in \cite{VANHOEIJ2000} and \cite{VANDERPUT2020}.

\bibliographystyle{plain}
\bibliography{DGT}

\end{document}